\documentclass%
 {ws-jnmp}

\usepackage{comments}
\newComments\DL{DL}{red}
\newComments\AL{AL}{blue}
\newComments\Ol{Ol}{green}

\usepackage{soul}
\usepackage{lscape}
\usepackage{mes}
\usepackage{graphicx}
\usepackage[all]{xy}

\usepackage{DLdef}

\begin{document}

\graphicspath{{./S3/}}


\markboth{S. Bouarroudj, P. Grozman, D. Leites} {Almost affine Lie
superalgebras}


\catchline{17}{Supp 01}{2010}{}{}

\copyrightauthor{Sofiane Bouarroudj, Pavel Grozman, Dimitry Leites}


\title{DEFINING RELATIONS OF ALMOST AFFINE (HYPERBOLIC)\\ LIE
SUPERALGEBRAS}


\author{SOFIANE BOUARROUDJ}
\address{Department of Mathematics, United Arab
Emirates University, Al Ain, PO. Box: 17551; United Arab Emirates\\
\email{Bouarroudj.sofiane@uaeu.ac.ae}}

\author{PAVEL GROZMAN}

\address{Equa
Simulation AB, Stockholm, Sweden\\
\email{pavel.grozman@bredband.net}}

\author{DIMITRY LEITES}

\address{Dept. of Math., Stockholm University, Roslagsv. 101,
Kr\"aftriket hus 6, SE-106 91,
 Stockholm, Sweden\\
\email{mleites@math.su.se}}

\maketitle

\begin{history}
\received{(20 February 2009)} 
\accepted{(20 May 2009)}
\end{history}

\begin{abstract}
For all almost affine (hyperbolic) Lie superalgebras, the defining
relations are computed in terms of their Chevalley generators.
\end{abstract}


\ccode{2000 Mathematics Subject Classification: 17B65}

\section{Introduction}

For the definition of Lie superalgebras with Cartan matrix, in
particular, of the almost affine Lie superalgebras, see \cite{CCLL}.
We recall the definition of Chevalley generators and the defining
relations expressed in terms of these generators, see \cite{GL2} and
a review \cite{BGLL} which also contains the modular case.

Here we list defining relations of almost affine Lie superalgebras
with indecomposable Cartan matrices classified in \cite{CCLL}. For
the Lie superalgebras whose Cartan matrices are symmetrizable and
without zeros on the main diagonal, the relations were known; they
are described by almost the same rules as for Lie algebras and they
are only of Serre type if the off-diagonal elements of such Cartan
matrices are non-positive.

Nothing was known about Lie superalgebras of the type we study in
this note: almost affine Lie superalgebras with indecomposable
Cartan matrices with zeros on the diagonal and without zeros but
non-symmetrizable.

Although presentation --- description in terms of generators and
relations --- is one of the accepted ways to represent a given
algebra, it seems that an {\it explicit} form of the presentation is
worth the trouble to  obtain only if this presentation is often in
need, or (which is usually the same) is sufficiently neat. The
Chevalley generators of simple finite dimensional Lie algebras over
$\Cee$ satisfy simple and neat relations (\lq\lq Serre relations'')
and are often needed for various calculations and theoretical
discussions. Relations between their analogs in the super case,
although not so neat (certain \lq\lq non-Serre relations'' appear),
are still tolerable, at least, for most Cartan matrices. One Lie
superalgebra may have several inequivalent Cartan matrices, and
although for some Cartan matrices of a given Lie superalgebra the
relations between Chevalley generators are rather complicated
(non-Serre ones), for some other Cartan matrices of the same algebra
there are only Serre relations, see \cite{GL2}; see \cite{RU} where
this fact is used.

What we usually need to know about defining relations is that there
are finitely many of them; hence the fact that some simple loop
superalgebras with Cartan matrix are {\bf not} finitely presentable
in terms of Chevalley generators was unexpected (although obvious as
an afterthought, see \cite{GL2}, where this was first published).

The non-Serre relations --- a bit too complicated to be used by
humans --- were of purely theoretic interest until recently
Grozman's package \texttt{SuperLie} (\cite{Gr}) made the task of
finding the explicit expression of the defining relations for many
types of Lie algebras and superalgebras a routine exercise for
anybody capable to use \texttt{Mathematica}. Unfortunately, such
calculations, performed up to certain degree in generators, are not
conclusive if not supported by cohomological arguments of the type
elucidated in \cite{Zus}: there might be some relations of higher
orders left unobserved.

There are several reasons for deriving the relations even if they
look awful:
\begin{enumerate}
  \item Computers easily swallow, often, what is non-palatable to humans (and
  the relation we derive in this note is precisely the case);
  \item even a rough description of the relations (the value of the highest
  degree, etc.) might be sufficient for deriving the next result
  based on this description
  (e.g., the description of the real forms of the Lie algebra
  $\fgl(\lambda)$ of \lq\lq matrices of complex size", see
  \cite{LSg} based on \cite{GL1}).
\end{enumerate}

Incompleteness of relations is irrelevant if it is clear that there
are infinitely many relations and their explicit form is not,
actually, needed, as in some cases (like $\fpsl(n|n)^{(1)}$)
considered in \cite{GL1}.

\section{Recapitulation}
Let the elements $X_i^\pm$, $d_i$,  and $h_i$, where $i=1,\dots,n$,
generate the Lie superalgebra denoted $\tilde \fg(A,
I)$,\index{$\tilde\fg(A, I):=\tilde \fg(A, I)$} where $I=(i_1, \dots
i_n)\in(\Zee/2)^n$ is a collection of parities ($p(e_j^\pm)=i_j$ and
$A=(A_{ij})_{i=1}^n$), see \cite{CCLL}.

\ssec{Serre relations, see [GL2]} Let $A$ be an $n\times n$ matrix.
We find the defining relations by induction on $n$ with the help of
the Hochschild--Serre spectral sequence (for its description for Lie
superalgebras, which has certain subtleties, see \cite{Po}). For the
basis of induction consider the following cases of Dynkin diagrams
with one vertex and no edges:
\begin{equation}
\label{3.4.1}
\begin{array}{ll}
 \mcirc \text{ or } \mbullet &\text{no relations, i.e., $\fg^{\pm}$
are free Lie superalgebras};\\
\mbullet &\ad_{X^{\pm}}^2(X^{\pm})=0;\\
{\motimes}& [X^{\pm}, X^{\pm}]=0.
\end{array}
\end{equation}
Set $\deg X_{i}^{\pm} = 0$ for $ 1\leq i\leq n-1$ and $\deg
X_{n}^{\pm} = \pm 1$. Let $\fg^{\pm} = \oplus \fg_{ i}^{\pm}$ and
$\fg= \oplus \fg_{i} $ be the corresponding $\Zee $-gradings. Set
$\fg_{\pm} =\fg^{\pm}/\fg_{0}^{\pm}$. From the Hochschild--Serre
spectral sequence for the pair $\fg_{ 0}^{\pm} \subset \fg ^{\pm}$
we get:
\begin{equation}
\label{3.4.2} H_{2}(\fg_{\pm})\subset H_{2}(\fg_{0}^{\pm})\oplus
H_{1}(\fg_{ 0}^{\pm}; H_{ 1}(\fg_{\pm}))\oplus H_{0}(\fg_ {
0}^{\pm}; H_{ 2}(\fg_{\pm})).
\end{equation}
It is clear that
\begin{equation}
\label{3.4.3} H_{1}(\fg_{\pm})= \fg_{1}^{\pm} , \;\;\; H_{
2}(\fg_{\pm}) = \wedge^{2}(\fg_ {1}^{\pm})/\fg_{ 2}^{\pm}.
\end{equation}
So, the second summand in~(\ref{3.4.2}) provides us with relations
of the form:
\begin{equation}
\label{3.4.4} \begin{array}{ll} (\ad_{X_{n}^{\pm}})^{k_{ni}}
(X_{i}^{\pm})=0&\text{if the $n$-th simple root is not}\;\;
{\motimes}\\
{} [X_{n}, X_{n}]=0&\text{if the $n$-th simple root is}
\;\;\motimes.
\end{array}
\end{equation}
while the third summand in (\ref{3.4.2}) is spanned by the
$\fg_{0}^{\pm}$-lowest vectors in
\begin{equation}
\label{3.4.6} \wedge^{2}(\fg_{1}^{\pm})/(\fg_{2}^{\pm} + \fg^{\pm}
\wedge^{2}(\fg_{1}^{\pm})).
\end{equation}

Let the matrix $B=(B_{kj})$ be given by the formula
\begin{equation}
\label{Boddrefl}B_{kj}=\begin{cases}
-\displaystyle\frac{2A_{kj}}{A_{kk}}& \text{~if~} A_{kk}\neq
0\text{~and~}
-\displaystyle\frac{2A_{kj}}{A_{kk}}\in \Zee_+, \\
1&\text{~if~}i_k=\od, A_{kk}=0,A_{kj}\neq 0,\\
0&\text{~if~}i_k=\od, A_{kk}=A_{kj}=0,\\
0&\text{~if~}i_k=\ev, A_{kk}=\ev,A_{kj}=0.\end{cases}
\end{equation}
The following proposition, its proof being straightforward,
illustrates the usefulness of our normalization of Cartan matrices
as compared with other options:

\sssbegin{Proposition} The numbers $k_{in}$ and $k_{ni}$ in
$(\ref{3.4.4})$ are expressed in terms of $B_{kj}$ as follows:
\begin{equation}
\label{srpm}
\begin{matrix} (\ad_{X_{k}^{\pm}})^{1+B_{kj}}(X_{j}
^{\pm})=0 &
\text{ for $k \neq j$} \\ & \\
\text{$[X_{i}^{\pm}, X_{i}^{\pm}]=0$} & \text{if $A_{ii}
=0$}.\end{matrix}
\end{equation}
\end{Proposition}
The relations (\ref{gArel_0}) and $(\ref{srpm})$ will be called {\it
Serre relations} of the Lie superalgebra $\fg(A)$.

\ssec{Non-Serre relations} These are relations that correspond to
the third summand in~(\ref{3.4.2}). Let us consider the simplest
case: $\fsl (m|n)$ in the realization with the system of simple
roots
\begin{equation}
\label{circ}
 \xymatrix@C=1em{
 \ffcirc\ar@{-}[r]&\dots\ar@{-}[r]&\ffcirc\ar@{-}[r]&
 \ffotimes\ar@{-}[r]&\ffcirc\ar@{-}[r]&\dots\ar@{-}[r]&
 \ffcirc
 }
\end{equation}

Then $H_2(\fg_{\pm})$ from the third summand in (\ref{3.4.2}) is
just $\wedge^2(\fg_{\pm})$. For simplicity, we confine ourselves to
the positive roots. Let $X_{1}$, \dots , $X_{m-1}$ and
$Y_{1}$,~\dots,~$Y_{n-1}$ be the root vectors corresponding to even
roots separated by the root vector $Z$ corresponding to the root
$\otimes$.

If $n=1$ or $m=1$, then $\wedge^2(\fg)$ is an irreducible $\fg_{\bar
0}$-module and there are no non-Serre relations. If $n\neq 1$ and
$m\neq 1$, then $\wedge^2(\fg)$ splits into 2 irreducible $\fg_{\bar
0}$-modules. The lowest component of one of them corresponds to the
relation $[Z, Z]=0$, the other one corresponds to the non-Serre-type
relation
\begin{equation}
\label{*} [[X_{m-1}, Z], [Y_{1}, Z]] =0. \end{equation}

If, instead of $\fsl (m|n)$, we would have considered the Lie
algebra $\fsl(m+n)$, the same argument would have led us to the two
relations, both of Serre type:
\[
\ad_Z^2(X_{m-1})=0, \qquad \ad_Z^2(Y_{1})=0.
\]

In what follows we give an explicit description of the defining
relations between the generators of the positive part of the almost
affine Lie superalgebras in terms of their Chevalley generators.

\section{Results (the proof is based on [Zus] and the induction as above)}
For $NS3_{33}$ (CM2), CM4): Serre relations only):
\begin{equation}\label{rels}
\renewcommand{\arraystretch}{1.4}
\begin{array}{ll}
CM1)&[x_1,\;x_1]=0,\quad [x_2,\;x_2]=0,\quad [x_3,\;x_3]=0,\\
& 2[[x_1,\;x_2],\;[x_2,\;[x_1,\;x_3]]]= a\,
[[x_1,\;x_2],\;[x_3,\;[x_1,\;x_2]]],\\
&
[[x_1,\;x_3],\;[x_2,\;[x_1,\;x_3]]]=2 a\, [[x_1,\;x_3],\;[x_3,\;[x_1,\;x_2]]],\\
& [[x_2\text{,$\, $}x_3],\;[[x_2,\;x_3],\;[x_2,\;[x_1,\;x_3]]]]=(3 a
+5)\, [[x_2,\;x_3],\;[[x_2,\;x_3],\;[x_3\text{,$\,
   $}[x_1,\;x_2]]]].\\
CM3)&[x_2,\;x_2]=0,\quad [x_1,\;[x_1,\;x_2]]=0,\quad
[x_1,\;[x_1,\;x_3]]=0,\quad [x_3,\;[x_2,\;x_3]]=0,\\
& [x_3,\;[x_3,\;[x_1,\;x_3]]]=
0,\\
& 10[[x_1,\;x_3],\;[[x_1,\;x_2],\;[x_2,\;x_3]]]=20(a + 2)\,
[[x_2,\;[x_1,\;x_3]],\;[x_3,\;[x_1\text{,$\,
   $}x_2]]]-\\
& (15a +13)\, [[x_3,\;[x_1,\;x_2]],\;[x_3\text{,$\,
   $}[x_1,\;x_2]]]- (15a+8)\, [[x_2,\;[x_1,\;x_3]]\text{,$\,
   $}[x_2,\;[x_1,\;x_3]]].\\
\end{array}
\end{equation}

For $NS3_{34}$ (CM2), CM4): Serre relations only):
\begin{equation}\label{rels2}
\renewcommand{\arraystretch}{1.4}
\begin{array}{ll}
CM1)&[x_1,\;x_1]=0,\quad [x_2,\;x_2]=0,\quad [x_3,\;x_3]=0,\\
& 2[[x_1,\;x_2],\;[x_2,\;[x_1,\;x_3]]]= a\,
[[x_1,\;x_2],\;[x_3,\;[x_1,\;x_2]]],\\
&
[[x_1,\;x_3],\;[x_2,\;[x_1,\;x_3]]]=2 a\, [[x_1,\;x_3],\;[x_3,\;[x_1,\;x_2]]]\\
CM3)&[x_2,\;x_2]=0,\quad [x_1,\;[x_1,\;x_2]]=0,\\
& [x_1,\;[x_1,\;x_3]]=0,\quad [x_3,\;[x_2,\;x_3]]=0,\quad
[x_3,\;[x_3,\;[x_3,\;[x_1,\;x_3]]]]=
   0,\\
 &
 14[[x_1,\;x_3],\;[[x_1,\;x_2],\;[x_2,\;x_3]]]=42(a +4)\,
   [[x_2,\;[x_1,\;x_3]],\;[x_3,\;[x_1\text{,$\,
   $}x_2]]]-\\
&2(14a +11)\, [[x_3,\;[x_1,\;x_2]],\;[x_3\text{,$\,
   $}[x_1,\;x_2]]]-(28a +15)\, [[x_2,\;[x_1,\;x_3]]\text{,$\,
   $}[x_2,\;[x_1,\;x_3]]].\\
\end{array}
\end{equation}

For $NS3_{35}$ (CM2), CM4): Serre relations only):
\begin{equation}\label{rels2}
\renewcommand{\arraystretch}{1.4}
\begin{array}{ll}
CM1)&[x_1,\;x_1]=0,\quad [x_2,\;x_2]=0,\quad [x_3,\;x_3]=0,\\
& 2[[x_1,\;x_2],\;[x_2,\;[x_1,\;x_3]]]= a\,
[[x_1,\;x_2],\;[x_3,\;[x_1,\;x_2]]]\\
&
[[x_1,\;x_3],\;[x_2,\;[x_1,\;x_3]]]=2 a\, [[x_1,\;x_3],\;[x_3,\;[x_1,\;x_2]]].\\
CM3)& [x_2,\;x_2]=0,\quad [x_1,\;[x_1,\;x_2]]=0,\quad
[x_1,\;[x_1,\;x_3]]=0,\quad [x_3,\;[x_2,\;x_3]]=0,\\
&
[x_3,\;[x_3,\;[x_3,\;[x_3,\;[x_1,\;x_3]]]]]= 0\\
& 6[[x_1,\;x_3],\;[[x_1,\;x_2],\;[x_2,\;x_3]]]=24(a +2)\,
[[x_2,\;[x_1,\;x_3]],\;[x_3,\;[x_1\text{,$\,
   $}x_2]]]-\\
&(15a+11)\, [[x_3,\;[x_1,\;x_2]],\;[x_3\text{,$\,
   $}[x_1,\;x_2]]]-(15a +8)\, [[x_2,\;[x_1,\;x_3]]\text{,$\,
   $}[x_2,\;[x_1,\;x_3]]].\\
   \end{array}
\end{equation}

For $NS3_{36}$ (CM2), CM3),CM4): Serre relations only):
\begin{equation}\label{rels2}
\renewcommand{\arraystretch}{1.4}
\begin{array}{ll}
CM1)&[x_1,\;x_1]=0,\quad [x_2,\;x_2]=0,\quad [x_3,\;x_3]=0,\\
& 3[[x_1,\;x_2],\;[[x_1,\;x_2],\;[x_2,\;[x_1,\;x_3]]]]= a\,
[[x_1,\;x_2],\;[[x_1,\;x_2],\;[x_3,\;[x_1,\;x_2]]]],
\\ &[[x_1,\;x_3],\;[[x_1,\;x_3],\;[x_2,\;[x_1,\;x_3]]]]= 3 a\,
[[x_1,\;x_3],\;[[x_1,\;x_3],\;[x_3,\;[x_1,\;x_2]]]],\\&
[[x_2,\;x_3],\;[x_2,\;[x_1,\;x_3]]]=(2a +5)\,
   [[x_2,\;x_3],\;[x_3,\;[x_1\text{,$\,
   $}x_2]]],\\
   &16[[x_2,\;x_3],\;[[x_1\text{,$\,
   $}x_2],\;[x_1,\;x_3]]]=-(15a +37)\,
   [[x_2,\;[x_1,\;x_3]],\;[x_2,\;[x_1\text{,$\,
   $}x_3]]]+\\
&(10a +14)\, [[x_2,\;[x_1,\;x_3]],\;
   [x_3,\;[x_1\text{,$\,
   $}x_2]]]-5(3a+1)\, [[x_3,\;[x_1,\;x_2]],\;[x_3\text{,$\,
   $}[x_1,\;x_2]]].\\
\end{array}
\end{equation}

For $NS3_{37}$ (CM2), CM3),CM4): Serre relations only):
\begin{equation}\label{rels2}
\renewcommand{\arraystretch}{1.4}
\begin{array}{ll}
CM1)& [x_1,\;x_1]=0,\quad [x_2,\;x_2]=0,\quad [x_3,\;x_3]=0,\\
& 3[[x_1,\;x_2],\;[[x_1,\;x_2],\;[x_2,\;[x_1,\;x_3]]]]= a\,
[[x_1,\;x_2],\;[[x_1,\;x_2],\;[x_3,\;[x_1,\;x_2]]]]\\
&[[x_1,\;x_3],\; [[x_1,\;x_3],\;[x_2,\;[x_1,\;x_3]]]]=
3 a\, [[x_1,\;x_3],\;[[x_1,\;x_3],\;[x_3,\;[x_1,\;x_2]]]]\\
& [[x_2,\;x_3],\;[[x_2,\;x_3],\;[x_2,\;[x_1,\;x_3]]]]=(3a +8)\,
[[x_2,\;x_3],\;[[x_2,\;x_3],\;[x_3\text{,$\,
   $}[x_1,\;x_2]]]]\\
\end{array}
\end{equation}

For the remaining almost affine Lie superalgebras, the defining
relations are only Serre ones regardless of symmetrizability of
their Cartan matrix.


\end{document}